\documentclass{article}
\usepackage[utf8x]{inputenc}

\usepackage{a4wide}
\usepackage{amscd}
\usepackage{amsmath}
\usepackage{amssymb}
\usepackage{amsthm}
\usepackage{array}
\usepackage{color}
\usepackage{epsfig}
\usepackage{graphics}
\usepackage{graphicx}
\usepackage[colorlinks=true,urlcolor=blue,linkcolor=blue,plainpages=false,pdfpagelabels
]{hyperref}
\usepackage{ifthen}
\usepackage{tabularx}
\usepackage{url}
\usepackage[all,cmtip]{xy}

\newcounter{uebno}

\newsavebox{\uebung}

\newcommand{\solutionsem}[1]{} 

\newtheorem{theorem}{Theorem}[section]

\newtheorem{definition}[theorem]{Definition}
\newtheorem{proposition}[theorem]{Proposition}
\newtheorem{corollary}[theorem]{Corollary}
\newtheorem{lemma}[theorem]{Lemma}
\newtheorem{fact}[theorem]{Remark}
\newtheorem{exemplu}[theorem]{Example}
\newtheorem{exercise}{Exercise}
\newtheorem{notation}[theorem]{Notation}

\newcommand{\bdfn}{\begin{definition}}
\newcommand{\edfn}{\end{definition}}
\newcommand{\bthm}{\begin{theorem}}
\newcommand{\ethm}{\end{theorem}}
\newcommand{\bprop}{\begin{proposition}}
\newcommand{\eprop}{\end{proposition}}
\newcommand{\bcor}{\begin{corollary}}
\newcommand{\ecor}{\end{corollary}}
\newcommand{\blem}{\begin{lemma}}
\newcommand{\elem}{\end{lemma}}
\newcommand{\bfact}{\begin{fact}}
\newcommand{\efact}{\end{fact}}
\newcommand{\bex}{\begin{exemplu}\begin{rm}}
\newcommand{\eex}{\end{rm}\end{exemplu}}
\newcommand{\bxc}{\begin{exercise}}
\newcommand{\exc}{\end{exercise}}
\newcommand{\bntn}{\begin{notation}}
\newcommand{\entn}{\end{notation}}

\newcommand{\be}{\begin{enumerate}}
\newcommand{\ee}{\end{enumerate}}
\newcommand{\bce}{\begin{center}}
\newcommand{\ece}{\end{center}}
\newcommand{\bi}{\begin{itemize}}
\newcommand{\ei}{\end{itemize}}
\newcommand{\bt}{\begin{tabular}}
\newcommand{\et}{\end{tabular}}
\newcommand{\beq}{\begin{equation}}
\newcommand{\eeq}{\end{equation}}
\newcommand{\ba}{\begin{array}} 
\newcommand{\ea}{\end{array}}
\newcommand {\bea} {\begin{eqnarray}}
\newcommand {\eea} {\end {eqnarray}}
\newcommand {\bua} {\begin{eqnarray*}}
\newcommand {\eua} {\end {eqnarray*}}

\newcommand{\se}{\subseteq}

\newcommand{\ds}{\displaystyle}

\newcounter{ct}

\def\R{{\mathbb R}}
\def\N{{\mathbb N}}

\def\R{{\mathbb R}}

\newcommand{\eps}{\varepsilon}

\newcommand{{\Fdo}}{F}

\title{Effective results on a fixed point algorithm for families of nonlinear mappings}
\author{Andrei Sipo\c s${}^{a,b}$\\[0.2cm]
\footnotesize ${}^a$Faculty of Mathematics and Computer Science, University of Bucharest,\\
\footnotesize Academiei 14, 010014 Bucharest, Romania\\[0.1cm]
\footnotesize ${}^b$Simion Stoilow Institute of Mathematics of the Romanian Academy,\\
\footnotesize P. O. Box 1-764, 014700 Bucharest, Romania\\[0.1cm]
\footnotesize E-mail: Andrei.Sipos@imar.ro
}
\date{}

\begin{document}

\maketitle

\begin{abstract}
We use proof mining techniques to obtain a uniform rate of asymptotic regularity for the instance of the parallel algorithm used by López-Acedo and Xu to find common fixed points of finite families of $k$-strict pseudocontractive self-mappings of convex subsets of Hilbert spaces. We show that these results are guaranteed by a number of logical metatheorems for classical and semi-intuitionistic systems.


\end{abstract}

\section{Introduction}

The class of $k$-strict pseudocontractions was introduced and studied by Browder and Petryshyn in \cite{BroPet67} in the context of Hilbert spaces. If $H$ is a Hilbert space, $C \se H$ is a convex subset and $k \in [0,1)$, then a mapping $T: C \to H$ is called a {\bf $k$-strict pseudocontraction} if for all $x,y \in C$ we have that:
\beq
\|Tx-Ty\|^2 \leq \|x-y\|^2 + k\|(x-Tx) - (y-Ty)\|^2.\label{k-ps}
\eeq
If we set $k:=0$ in the above inequality, we obtain the condition $\|Tx-Ty\| \leq \|x-y\|$, which states that the mapping $T$ is nonexpansive. The search of algorithms for finding fixed points of nonexpansive self-mappings of subsets of metric spaces belonging to various established classes has been a longstanding research program in the field of nonlinear analysis.

These algorithms are usually iterative in their nature -- one begins from an arbitrary starting point and proceeds to repeatingly apply to it a transformation which is derived from the operator $T$. In this way, the algorithm produces a sequence $(x_n)_{n \in \N} \se H$. Frequently, the first major intermediate step of the convergence proof is a result of {\it asymptotic regularity} -- that is, a statement of the form
$$\lim_{n \to \infty} \|x_n - Tx_n\| =0.$$
In this case, we also say that $(x_n)_n$ is an {\it approximate fixed point sequence} for $T$. Afterwards, using demiclosedness principles and/or various hypotheses of compactness, followed by considerations regarding the geometry of the space in question, one may obtain that the sequence $(x_n)_{n \in \N} \se H$ converges, weakly or strongly, to a fixed point of $T$.

Proof mining is a research program introduced by U. Kohlenbach in the 1990s (\cite{Koh08} is a comprehensive reference), which aims to obtain explicit quantitative information (witnesses and bounds) from proofs of an apparently ineffective nature. This paradigm in applied logic has successfully led so far to obtaining some previously unknown effective bounds, primarily in nonlinear analysis and ergodic theory. A large number of these are guaranteed to exist by a series of logical metatheorems which cover general classes of bounded or unbounded metric structures -- see \cite{Koh05,GerKoh06,GerKoh08}. As an example of such a piece of quantitative information, let us define a {\it rate of convergence} for a real-valued sequence $(a_n)_{n \in \N}$ that has the limit $a \in \R$ to be a function $\Theta : (0, \infty) \to \N$ such that for all $\eps >0$ and all $n \geq \Theta(\eps)$, we have that $\|a_n -a \|\leq \eps$. If the real-valued sequence is the sequence $(\|x_n - Tx_n\|)$ introduced before, then $\Theta$ is called a {\it rate of $T$-asymptotic regularity} for the sequence $(x_n)$. General metatheorems of proof mining usually guarantee that from the methods of proof used in nonlinear analysis one can obtain (``extract'') computable rates of asymptotic regularity for the usual iterative schemas -- see, for example, \cite{Leu07,Leu14}. In particular, one such recent result of Ivan and Leu\c stean \cite{IvaLeu15}, analysing a proof by Marino and Xu \cite{MarXu07}, states that for the Krasnoselski iteration of a single $k$-strict pseudocontractive self-mapping $T$ of a convex subset of a Hilbert space, there exists such a computable rate of $T$-asymptotic regularity that is quadratic in $\frac{1}{\eps}$.

That being said, such an area of investigation closely related to the kind of iterative algorithms mentioned before has been the problem of finding a common fixed point of a (finite or infinite) family $(T_i)_i$ of self-mappings of a subset $C$ like above. An iterative scheme that is useful in the case of a finite family $(T_i)_{1 \leq i \leq N}$ is the {\it parallel algorithm}, defined as follows. Let $x$ be in $C$ and $(t_n)_{n \in \N} \subseteq (0,1)$. For each $i \in \{1,\ldots,N\}$, let $(\lambda_i^{(n)})_{n\in\N}$ be a sequence of positive real numbers such that, for any $n \in \N$:
$$\sum_{i=1}^N \lambda_i^{(n)}=1.$$
Write, for all $n \geq 0$:
$$A_n:=\sum_{i=1}^N \lambda_i^{(n)} T_i.$$
Let $(x_n)_{n\in\N}$ be the sequence defined by:
$$x_0:=x$$
$$x_{n+1}:=t_nx_n + (1-t_n) A_n x_n$$
Then the sequence $(x_n)_{n\in\N}$ is the output of the parallel algorithm associated with the inputs $T$, $x$, $(t_n)$ and $(\lambda_i^{(n)})$.

Two remarks are in order here. Firstly, we see that the case $N=1$ represents the well-known Mann iteration for finding a fixed point of a self-mapping and therefore, all the results pertaining to this algorithm immediately transfer to the case of a single mapping (i.e. the one treated in \cite{IvaLeu15}). Secondly, we note that there exists another (equivalent) convention when working with Mann-like algorithms, pairing the $t_n$ with the application of the appropriate mapping, i.e. the formula above would be:
$$x_{n+1}:=t_nA_nx_n + (1-t_n) x_n$$
We use the ``$t_nx_n$'' convention, in the description of the parallel algorithm and also further below, when formalizing the passage from nonexpansive to strictly pseudocontractive mappings, as it is the one used in \cite{MarXu07,LopXu07}. One should be careful to check the convention used when comparing different hypotheses and convergence results.

When considering algorithms for finite families such as the one above, the intermediate result obtained during the proof of the convergence theorem will be still one of ``asymptotic regularity'', though one pertaining to the map(s) constructed as a byproduct of the algorithm (here, the $A_n$'s), i.e. that:
$$\lim_{n \to \infty} \|x_n - A_nx_n\| =0.$$
Given that $A_n$ varies with $n$, such a result does not mean {\it a priori} that $(x_n)_n$ is an approximate fixed point sequence for any mapping -- certainly not one given by the problem data. Therefore, what is actually relevant to the proof mining program is an asymptotic regularity related to the relevant mappings of the problem -- that is, the $T_i$'s. One might look for an associated rate of convergence for the statements:
$$\lim_{n \to \infty} \|x_n - T_ix_n\| =0,$$
for each $i \in \{1,...,N\}$. A concrete extraction of such a rate can be found, for example, in \cite{KhaKoh13}, for the Kuhfittig iteration.

López-Acedo and Xu, in 2007, have found sufficient conditions so that the parallel algorithm weakly converges to a fixed point of a finite family of strictly pseudocontractive mappings. Their result, using the notations introduced above, is expressed as follows.

\bthm[López-Acedo \& Xu (2007), {\cite[Theorem~3.3]{LopXu07}}]\label{lopxu}
Let $H$ be a Hilbert space and $C\se H$ a closed, convex set. Let $N \geq 1$, $(T_i : C \to C)_{1 \leq i \leq N}$ a family of mappings and $(k_i)_{1 \leq i \leq N} \se (0,1)$ such that each $T_i$ is a $k_i$-strict pseudocontraction. Suppose that $\bigcap_{i=1}^N Fix(T_i) \neq\emptyset$. Set $k:=\max_{1\leq i \leq N} k_i$. Let $x$ be in $C$, $(t_n)_{n \in \N} \subseteq [k,1]$ be such that
$$\sum_{n=0}^\infty (t_n-k)(1-t_n) = \infty.$$
Impose the conditions
$$\inf_{i,n} \lambda_i^{(n)} >0$$
and
$$\sum_{j=0}^\infty \sqrt{\sum_{i=1}^N |\lambda_i^{(j+1)}-\lambda_i^{(j)}|} < \infty$$
on $(\lambda_i^{(n)})$.
Then the the parallel algorithm associated with the inputs $T$, $x$, $(t_n)$ and $(\lambda_i^{(n)})$ weakly converges to a common fixed point of the family $(T_i : C \to C)_{1 \leq i \leq N}$.
\ethm

Our goal will be to obtain rates of asymptotic regularity for this instance of the parallel algorithm -- that is, a rate of convergence for each sequence $(\|x_n-T_ix_n\|)_{n \in \N}$, with the sequence $(x_n)_{n\in\N}$ being defined as before.

\section{A discussion of the proof from a logical standpoint}

\subsection{The logical systems and metatheorems involved}

The ultimate grounds in extracting quantitative information from a proof are the general logical metatheorems of proof mining, developed by Kohlenbach \cite{Koh05} and by Gerhardy and Kohlenbach \cite{GerKoh06,GerKoh08}. The hypotheses of the metatheorems that usually play the crucial r\^ole in determining whether the extraction can be carried out for the case at hand are firstly, the possibility of the formalization of the proof in one of few studied logical systems and secondly, a suitable form of low logical complexity for the conclusion of the theorem. In the sequel, we shall present the logical prerequisites which will allow for a rough presentation and understanding of the metatheorems which are relevant to our results.

Consider the system $\mathcal{A}^\omega$ of weakly extensional classical analysis in all finite types -- that is, the theory weakly extensional classical (Peano) arithmetic in all finite types supplemented by the quantifier-free axiom schema of choice and the axiom schema of dependent choice. This system has a power comparative to the first-order, two-sorted theory usually denoted by $\text{\bf Z}_\text{\bf 2}$ or ``full second-order arithmetic''. For a detailed exposition of the system $\mathcal{A}^\omega$, including the way it represents real numbers, one is invited to consult \cite{Koh08}. We note that in this case, the real numbers are represented by the type $1=0(0)$ and the relations $=_\R$ and $\leq_\R$ are implemented by purely universal formulas and $<_\R$ by a purely existential formula. As many ordinary definitions in analysis can freely interchange ``$<$'' with ``$\leq$'' (take, for example, the usual definition of a Cauchy sequence), this universal/existential duality will help us use formulas of lower logical complexity when the framework demands it.

Now, in order to be able to express properties of Hilbert spaces, we have to extend this logical system. We shall closely follow \cite{Koh08}.

The first such extension, $\mathcal{A}^\omega[X,\|\cdot\|]$, allows us to speak about normed spaces. The set of types for this system will be generated by two ``primitive'' types, the type $0$ of natural numbers and a new abstract type $X$, representing elements from our space, forming a free algebra with a single binary operation $\to$, representing function types. (We will write $\tau(\rho)$ for $\rho\to\tau$.) For such a type $\rho$, we define the type $\hat{\rho}$ by replacing all occurences of $X$ in $\rho$ by $0$. Here we also briefly define the notion of {\it degree} for such a type, as follows. We say that $\rho$ has degree:
\bi
\item $\leq 1$ if $\rho = 0$ or $\rho = 0(0)...(0)$;
\item $1^*$ if $\hat{\rho}$ has degree $\leq 1$;
\item $(0,X)$ if $\rho=X$ or $\rho = X(0)...(0)$;
\item $(1,X)$ if $\rho=X$ or $\rho = X(\rho_n)...(\rho_1)$ where each $\rho_i$ has degree $\leq 1$ or $(0,X)$.
\ei

Also, we add new constants for the various operations common to normed spaces, i.e. $0_X$ and $1_X$ of type $X$, $+_X$ of type $X(X)(X)$, $-_X$ of type $X(X)$, $\cdot_X$ of type $X(X)(1)$ (remember that $1=0(0)$ is the type of real numbers) and $\|\cdot\|_X$ of type $1(X)$. We allow infix notation and the ``syntactic sugar'' of writing $x -_X y$ for $x +_X (-_X y)$. Finally, we add the following axioms:
\be
\item the equational, and hence purely universal, axioms for vector spaces;
\item $\forall x^X (\|x-_X x\|_X =_\R 0_\R )$;
\item $\forall x^X y^X (\|x -_X y\|_X =_\R \|y -_X x\|_X)$;
\item $\forall x^Xy^Xz^X (\|x -_X z\|_X \leq_\R \|x-_X y\|_X +_\R \|y-_Xz\|_X)$;
\item $\forall \alpha^1 x^X y^X (\|\alpha x-_X\alpha y\|_X =_\R \|\alpha\|_\R \cdot_\R \|x-_X y\|_X$;
\item $\forall \alpha^1 \beta^1 x^X (\|\alpha x -_X \beta x\|_X =_\R |\alpha -_\R \beta |_\R \cdot_\R \|x\|_X$;
\item $\forall x^X \forall y^X \forall u^X \forall v^X (\|(x+_X y) -_X (u+_X v)\|_X \leq_\R \|x -_X u\|_X +_\R \|y-_X v\|_X)$;
\item $\forall x^X y^X (\|(-_X x) -_X (-_X y)\|_X =_\R \|x-_X y\|_X)$;
\item $\forall x^X y^X (|\|x\|_X -_\R \|y\|_X|_\R \leq_\R \|x-_X y\|_X)$;
\item $\|1_X\|_X =_\R 1_\R$.
\ee

Note that the equality relation $x^X =_X y^X$ which is necessarily used in the expression of the vector space axioms is syntactically defined as $\|x -_X y\|_X =_\R 0_\R$. We define the equality for higher types as in the system $\mathcal{A}^\omega$, as extensional equality reducible to $=_0$ and $=_X$.

An issue when adding new constant symbols is their extensionality -- roughly, as the base system admits only a quantifier-free rule of extensionality, it is not {\it prima facie} clear that for a new function symbol $f$ one can prove in the new system a statement of the form
$$\forall x_1...\forall x_n\forall y_1...\forall y_n (\bigwedge_i x_i = y_i \to f(x_1,...,x_n) = f(y_1,...,y_n))$$
Some axioms above, like the eighth one, are written in this way purely to minimize the effort in writing such an extensionality proof; the rest of them yield it more readily in their classical forms. The result is that all new function symbols are provably extensional. The last axiom is added solely to ensure the non-triviality of the formalized space.

The theory of Hilbert spaces (actually, pre-Hilbert spaces or inner product spaces, since we do not consider completeness here) together with a nonempty unbounded convex subset, $\mathcal{A}^\omega[X,\langle \cdot,\cdot \rangle, C]_{-b}$ is obtained from $\mathcal{A}^\omega[X,\|\cdot\|]$ by adding the new constants $c_X$ of type $X$ and $\chi_C$ of type $0(X)$ and the following axioms (where $\tilde{\alpha}$ is $0$ if $\alpha<0$ and $1$ if $\alpha>1$, otherwise it is the same as $\alpha$; it follows from the representation of real numbers that this construction can be effectively coded):
\be
\item $\forall x^X y^X (\|x +_X y\|_X^2 +_\R \|x-_X y\|_X^2 =_\R 2_\R \cdot_\R (\|x\|_X^2 +_\R \|y\|_X^2)$ (the parallelogram law, which is a necessary and sufficient condition for a normed structure to arise from an inner product structure, the last of which being necessarily unique);
\item $\forall x^X y^X \alpha^1 (\chi_C(x) =_0 0 \land \chi_C(y) =_0 0 \to \chi_C((1-_\R \tilde{\alpha}) \cdot_X x +_X \tilde{\alpha} \cdot_X y) =_0 0)$;
\item $\chi_C(c_X) =_0 0$;
\item $\forall x^X (\chi_C(x) \leq_0 1)$.
\ee

Now, if $(H,\langle\cdot,\cdot\rangle)$ is an inner product space and $C \se X$ is a nonempty convex subset, we can associate a set-theoretic model $\mathcal{S}^{\omega,X}$ in all finite types to this data by putting $\mathcal{S}_0 := \N$, $\mathcal{S}_X := H$ and $\mathcal{S}_{\tau(\rho)} := \mathcal{S}_{\tau}^{\mathcal{S}_\rho}$, assigning to any language constant its standard value, except for $1_X$ and $c_X$ who take arbitrary values that are of norm $1$ and are in $C$, respectively. We say that a sentence of our logical language is modeled by such a triple $(H,\langle\cdot,\cdot\rangle,C)$ iff it is satisfied in the usual Tarskian sense by all the possible models associated to it (i.e., we let $1_X$ and $c_X$ vary arbitrarily within their limits).

We say that a formula in our language is a $\forall$-formula (resp. an $\exists$-formula) iff it is formed by adjoining a list of universal (resp. existential) quantifiers over variables of types of degree $1^*$ or $(1,X)$ to a quantifier-free formula.

The following theorem is an adaptation of one by Gerhardy and Kohlenbach (see \cite{GerKoh08} and \cite[Corollary 17.71]{Koh08}).

\bthm
Let $P$ be $\N$, $\N^\N$ or $\N^{\N \times \N}$, $K$ an $\mathcal{A}^\omega$-definable compact metric space and $\rho$ a type of degree $1^*$, $B_\forall(u,y,z,f,n)$ a $\forall$-formula and $C_\exists(u,y,z,f,N)$ an $\exists$-formula (we assume that the free variables of the two formulas are at most the variables written as their arguments, such that the formula below shall be a sentence). Suppose that $\mathcal{A}^\omega[X,\langle \cdot,\cdot \rangle, C]_{-b}$ proves that:
$$\forall u \in P \forall M \in \N \forall k \in [0, 1-\frac{1}{M+1}] \forall y \in K \forall z^C \forall f^{C \to C} (\text{$f$ is $k$-strictly psc.}\land\forall n\in\N B_\forall \to \exists N \in \N C_\exists).$$
Then there exists a computable functional $\Phi : P \times \N \times \N \to  \N$ such that in all models $\mathcal{S}^{\omega, X}$ such that for all $x \in P$, $b \in \N$, $M \in \N$, $y\in K$, $z\in C$ and $f:C \to C$ a $\left(1-\frac{1}{N+1}\right)$-strictly pseudocontractive mapping, one has that:
$$\|z\| \leq b \land \|z-f(z)\| \leq b \land \forall n \leq \Phi(x,M,b) B_\forall \to \exists N \leq \Phi(x,M,b) C_\exists.$$
\ethm

\begin{proof}
One must ensure that $f$ is majorizable and that the additional premise of $k$-strict pseudocontractivity is universal. The latter is immediate, and the former follows because a $k$-strict pseudocontraction is Lipschitz of constant $\frac{1+k}{1-k}$ (see \cite{MarXu07} for the proof). Another byproduct of our hypotheses is that, by the remarks in \cite[p. 394]{Koh08}, $f$ is provably extensional and we may use its extensionality ``axiom'', similarly to the functions that belong to the definition of an inner product space.
\end{proof}

\bfact
Instead of single premises with single universal quantifiers, one might also have finite conjunctions of such. Also, one might replace the universally quantified variables at the beginning with finite tuples of them. (By this argument, the $k$, being taken from a compact definable space, is assimilated with the $y$ and therefore disappears as a separate dependency of the bound.)
\efact

The above remark will be of use in our case, since we will be dealing with several self-mappings of the convex set $C$.

Gerhardy and Kohlenbach \cite{GerKoh06} have also developed a series of metatheorems for so-called ``semi-intuitionistic'' systems -- that is, extensions of intuitionistic (``Heyting'') arithmetic by several principles that can be characterized as non-constructive. The analogous system to $\mathcal{A}^\omega[X,\langle \cdot,\cdot \rangle, C]_{-b}$ will be denoted by $\mathcal{A}_i^\omega[X,\langle \cdot,\cdot \rangle, C]_{-b}$ and will be based, instead, on fully extensional Heyting arithmetic in all finite types together with the full axiom of choice. The rest of the construction is analogous. We might consider, in addition, the schema of comprehension for negated formulas:
$$CA_\neg: \quad (\exists \Phi \leq \lambda x.1)(\forall y)(\Phi(y) =_0 0 \leftrightarrow \neg A(y))$$
in order to obtain the following metatheorem:

\bthm
Let $P$ be $\N$, $\N^\N$ or $\N^{\N \times \N}$, $K$ an $\mathcal{A}_i^\omega$-definable compact metric space and $B(u,y,z,f)$ and $C(u,y,z,f,N)$ be arbitrary formulas (using the same convention about free variables as before). Suppose that $\mathcal{A}_i^\omega[X,\langle \cdot,\cdot \rangle, C]_{-b} + CA_\neg$ proves that:
$$\forall u \in P \forall M \in \N \forall k \in [0, 1-\frac{1}{M+1}]  \forall y \in K \forall z^C \forall f^{C \to C} (\text{$f$ is $k$-strictly psc.}\land\neg B \to \exists N \in \N C).$$
Then there exists a computable functional $\Phi : P \times \N \times \N \to  \N$ such that in all models $\mathcal{S}^{\omega, X}$ such that for all $x \in P$, $b \in \N$, $M \in \N$, $y\in K$, $z\in C$ and $f:C \to C$ a $\left(1-\frac{1}{M+1}\right)$-strictly pseudocontractive mapping, one has that:
$$\exists N \leq \Phi(x,M,b) ( \neg B \to  C).$$
\ethm

\subsection{An analysis of the proof}

The last metatheorem of the previous section is proven using the modified realizability interpretation, a proof interpretation devised by G. Kreisel \cite{Kre59}. Modified realizability interprets intuitionistic logic, and hence cannot handle proofs which use the principle of excluded middle. These proof interpretations are translations from a given system to a usually quantifier-free system (but which for convenience is embedded into a ``full'' one) which also take notice, along the Hilbert-style proof tree, of ``realizers'' -- terms which can be fit into variables of the translated formula in order to witness it and to provide additional information about it. We can then transfer this information backwards in order to obtain realizers and bounds for the original formula, as it can be seen in the metatheorems above. It can be easily glimpsed that a soundness proof for the translation must proceed by way of an induction along the proof rules, providing realizers for axioms and ways of composing existing realizers for rules such as the modus ponens rule. A peculiar feature of these ``interpretations'' is that sometimes additional and well-known principles have a trivial translation that admits a similarly trivial term (like an identity function) as a valid realizer. We say that such a principle is {\it interpreted} by the given interpretation; we may also say that it is ``interpreted by itself''. Examples include the choice and comprehension principles stated in the last section.

Another such proof interpretation is G\"odel's functional (or ``Dialectica'') interpretation, together with its bar-recursion extension by Spector \cite{Spe62}, which also interprets only intuitionistic logic, with the difference that it can be more readily combined with the double negation translation from classical logic or arithmetic to their intuitionistic variants. The reason that this combination works in that Dialectica interprets the Markov principle:
$$\neg\neg\exists x A_0 (x) \to \exists x A_0 (x),$$
where $A_0$ is a quantifier-free formula that may have other free variables besides $x$. Such composition of interpretations is used to prove the first metatheorem. Not surprisingly, there is a trade-off in this approach, namely that, unlike modified realizability which can generally handle all sort of logical formulas, this one can only interpret those at roughly the $\forall\exists$ level, which can be seen in the restrictions imposed on the formulas in the first metatheorem of the previous section. Both kinds of interpretation are used in their monotone variants, introduced by Kohlenbach in \cite{Koh96,Koh98}.

Now, the first result of asymptotic regularity that we seek to quantify is the following:

\bthm
Using the notations from the Introduction, assume that $\sum_{n=0}^\infty (t_n-k)(1-t_n)=\infty$ and that $\sum_{j=0}^\infty \sum_{i=1}^N |\lambda_i^{(j+1)}-\lambda_i^{(j)}|<\infty$. Then $\lim_{n \to \infty} \|x_n -A_nx_n\|=0$.
\ethm

The statement above is a convergence result, which can be written in the standard expanded logical form of:
$$\forall \eps \exists N_\eps \forall n\geq N_\eps \|x_n -A_nx_n\| \leq \eps,$$
which is clearly seen to be a $\forall\exists\forall$ statement. We may call such a result {\it quantified} when we have found a computable functional $\Phi$ that gives the $N_\eps$ in terms of $\eps$ -- and such a functional is called a {\it rate of convergence} for the sequence. (We briefly remark that for a convergence result, a bound for $N_\eps$ would also be a realizer.)

Such a statement, however, is suitable only for the metatheorem which uses modified realizability. However, if one examines the proof of the above result from \cite{LopXu07}, one can see that it makes some use of the law of excluded middle, and hence the use of the first metatheorem is also not possible. Still, we shall find out that a rate of convergence can actually be extracted from the proof. This happens because of the way the proof is actually built. The first part of the proof, the one which uses excluded middle, actually establishes the following result:

\blem
With the same notations, we have that $\liminf_{n \to \infty} \|x_n -A_nx_n\|=0$.
\elem

Because the sequence in question is composed of real positive numbers, $0$ being a liminf of it is actually equivalent to $0$ being a limit point of the sequence, which can be written as:
$$\forall \eps \forall n \exists N\geq n \|x_N -A_Nx_N\| < \eps,$$
i.e. in a $\forall\exists$ form (note that we have used ``$<$'' instead of ``$\leq$'', following the previous discussion about the representation of real numbers).

One can, therefore, using the functional interpretation metatheorem, find a bound on the $N$ in the sentence above, which will turn it into a monotonely universal sentence. Such a sentence can be given as a hypothesis for the remained of the proof, which uses no non-constructive principles, and hence can be analysed using modified realizability. We note that a similar problem has previously been encountered and solved by Leu\c stean \cite{Leu14} for the Ishikawa iteration. (There is a hope that such cases may be more smoothly analysed using Oliva's hybrid interpretation \cite{Oli12}, which allows one to use in the same proof modified realizability and functional interpretation combined with negative translation in a modular way.)

All that remains now is to see that the metatheorems can be properly applied. By the remark in the previous section, one can replace the universal premise with a finite conjunction, i.e. add universal axioms to the framework. However, the premises that are at work here involve the convergence or divergence of series involving the families $(t_n)_n$ and $(\lambda_i^{(n)})_{i,n}$, which are clearly not universal. One solves this problem by adding rates of convergence and/or divergence to the statements and so turning them into universal ones. That is, we will take $\theta:\N\to\N$ be a rate of divergence for the series
$$\sum_{n=0}^\infty (t_n-k)(1-t_n),$$
i.e., for all $N \in \N$,
$$\sum_{n=0}^{\theta(N)} (t_n-k)(1-t_n) \geq N.$$
and $\gamma:(0,\infty)\to\N$ be a Cauchy modulus for the series
$$\sum_{j=0}^\infty \sum_{i=1}^N |\lambda_i^{(j+1)}-\lambda_i^{(j)}|,$$
i.e., for all $n \in \N$ and $\eps>0$,
$$\sum_{j=\gamma(\eps)+1}^{\gamma(\eps)+n} \sum_{i=1}^N |\lambda_i^{(j+1)}-\lambda_i^{(j)}| \leq \eps.$$
One gets in this way two new potential parameters for the eventual rate of convergence. However, we shall not be flooded with such parameters, as $(t_n)_n$ and $(\lambda_i^{(n)})_{i,n}$ are drawn out from compact intervals ($[k,1]$ and $[0,1]$, respectively) and so they will not take part in the final expression. The only remnants of them will be the $\theta$ and the $\gamma$.

We now, as announced, use the obtained bound on the $N$ (the ``modulus of liminf'') as an additional universal premise for the remainder of the proof, which is purely constructive, to get the desired rate of convergence for the asymptotic regularity sequence.

There is one final issue to discuss here: the extraction of a rate of convergence for the statement
$$\lim_{n \to \infty} \|x_n - T_ix_n\| =0$$
from the one which regards the statement
$$\lim_{n \to \infty} \|x_n - A_nx_n\| =0.$$

To do such an extraction, one might use the result which is implicit in the whole algorithm, namely that a fixed point of such a convex combination $A$ of mappings is a fixed point of each mapping $T_i$. Fix an $i$. Then the statement just expressed can be written as:
$$\forall q \in C (\forall \delta>0 \|q-Aq\| \leq \delta \to \forall \eps >0 \|q-T_iq\| < \eps)$$
or as:
$$\forall q \in C \forall \eps>0 \exists \delta>0( \|q-Aq\| \leq \delta \to  \|q-T_iq\| < \eps)$$
We can apply, again, on this $\forall\exists$ statement, the general metatheorem of proof mining in order to get the $\delta$ as a function of $\eps$. Then the rate of $T_i$-asymptotic regularity is immediately obtained.

\section{The main results}

\subsection{The case of nonexpansive mappings}

From now on, we will use the established notations and hypotheses from the definition of the parallel algorithm and from Theorem~\ref{lopxu}. Also, we will assume that there is a $b>0$ and a $p \in \bigcap_{i=1}^N Fix(T_i)$ such that $\|x\|\leq b$ and $\|x-p\|\leq b$. Let $\theta:\N\to\N$ and $\gamma:(0,\infty)\to\N$ be defined as in the last section.
Let $a>0$ be such that
$$a \leq \inf_{i,n} \lambda_i^{(n)}.$$

We will do the proof for the case $k=0$ at first (i.e. for nonexpansive mappings) and then, using a lemma of Browder and Petryshyn \cite{BroPet67}, derive the general result for strict pseudocontractions.

\blem\label{l22}
For all $n$, $t_n(1-t_n)\|x_n - A_nx_n\|^2 \leq \|x_n-p\|^2 - \|x_{n+1}-p\|^2$.
\elem

\begin{proof}
See \cite[(3.9)]{LopXu07}, for $k:=0$.
\end{proof}

\blem\label{l2}
The sequence $(\|x_n-p\|)_{n \in \N}$ is nonincreasing and for all $n$, $\|x_n-p\| \leq b$.
\elem

\begin{proof}
Immediate consequence of Lemma~\ref{l22} and of the fact that $(t_n)_{n \in \N} \subseteq (0,1)$.
\end{proof}

\bprop[``modulus of liminf'']\label{liminf}
Set, for all $\eps>0$ and $m \in \N$:
$$\Delta_{b,\theta}(\eps,m):=\theta(m+\left\lceil\frac{b^2}{\eps^2}\right\rceil).$$
Then, for all $\eps>0$ and $m \in \N$, there is an $N \in [m, \Delta_{b,\theta}(\eps,m)]$ such that $\|x_N - A_Nx_N\| \leq \eps$.
\eprop

\begin{proof}
Set $\ds \Psi:=\sum_{n=m}^{\Delta_{b,\theta}(\eps,m)} t_n(1-t_n)\|x_n - A_nx_n\|^2$. We have, using (i), that:
\begin{align*}
\Psi &\leq \sum_{n=0}^{\Delta_{b,\theta}(\eps,m)} t_n(1-t_n)\|x_n - A_nx_n\|^2 \\
&\leq \sum_{n=0}^{\Delta_{b,\theta}(\eps,m)} (\|x_n-p\|^2 - \|x_{n+1}-p\|^2) \\
&= \|x_0 - p\|^2 - \|x_{\Delta_{b,\theta}(\eps,m)+1}-p\|^2 \\
&\leq \|x - p\|^2 \\
&\leq b^2.
\end{align*}
We argue by contradiction. Suppose, then, that for all $n \in [m, \Delta_{b,\theta}(\eps,m)]$, $\|x_n - A_nx_n\| > \eps$. Then:
\begin{align*}
\Psi &> \eps^2 \sum_{n=m}^{\Delta_{b,\theta}(\eps,m)} t_n(1-t_n) \\
&= \eps^2 (\sum_{n=0}^{\Delta_{b,\theta}(\eps,m)} t_n(1-t_n) - \sum_{n=0}^{m-1} t_n(1-t_n))\\
&\geq \eps^2 (\sum_{n=0}^{\theta(m+\lceil\frac{b^2}{\eps^2}\rceil)} t_n(1-t_n) - m)\\
&\geq \eps^2 (m+\lceil\frac{b^2}{\eps^2}\rceil - m)\\
&= \eps^2 \cdot\lceil\frac{b^2}{\eps^2}\rceil \\
&\geq b^2.
\end{align*}
Since we have by now proven that $\Psi \leq b^2$ and $\Psi > b^2$, we have obtained the desired contradiction.
\end{proof}

\blem
Each $A_n$ is a nonexpansive mapping.
\elem

\begin{proof}
We need to prove that a convex combination of nonexpansive mappings is nonexpansive. It is enough to prove this for two mappings, as the general case follows by induction. Let $T_1$ and $T_2$ be the two mappings and $T:= (1-t)T_1 + tT_2$ be their convex combination. But if $x$ and $y$ are two points, then:
\begin{align*}
\|Tx-Ty\| &= \|(1-t)(T_1x-T_1y) + t(T_2x-T_2y)\| \\
&\leq (1-t)\|T_1x-T_1y\| + t\|T_2x-T_2y\| \\
&\leq (1-t)\|x-y\| + t\|x-y\| \\
&= \|x-y\|.
\end{align*}
\end{proof}

\blem\label{lc}
We have that:
\be
\item $\|x_n-A_{n+1}x_{n+1}\| \leq (1-t_n)\|x_n-A_nx_n\| + \|x_{n+1}-A_{n+1}x_{n+1}\|$;
\item $\|A_nx_n - A_{n+1}x_{n+1}\| \leq (1-t_n)\|x_n - A_nx_n\| + \|A_{n+1}x_{n+1} - A_nx_{n+1}\|$;
\item $\|x_{n+1}-A_{n+1}x_{n+1}\| \leq (1-t_n)\|x_n-A_nx_n\| + t_n\|x_{n+1}-A_{n+1}x_{n+1}\| + (1-t_n)\|A_{n+1}x_{n+1} - A_nx_{n+1}\|$;
\item $\|x_{n+1}-A_{n+1}x_{n+1}\| \leq \|x_n-A_nx_n\| + \|A_{n+1}x_{n+1} - A_nx_{n+1}\|$;
\ee
\elem
\begin{proof}
\be
\item We have that:
\begin{align*}
\|x_n-A_{n+1}x_{n+1}\| &\leq \|x_n-x_{n+1}\| + \|x_{n+1}-A_{n+1}x_{n+1}\|\\
&= (1-t_n)\|x_n-A_nx_n\| + \|x_{n+1}-A_{n+1}x_{n+1}\|.
\end{align*}
\item We have that:
\begin{align*}
\|A_nx_n - A_{n+1}x_{n+1}\| &= \|A_nx_n - A_nx_{n+1} - (A_{n+1}x_{n+1} - A_nx_{n+1})\|\\
&\leq \|A_nx_n - A_nx_{n+1}\| + \|A_{n+1}x_{n+1} - A_nx_{n+1}\|\\
&\leq  \|x_n-x_{n+1}\| + \|A_{n+1}x_{n+1} - A_nx_{n+1}\|\\
&= (1-t_n)\|x_n - A_nx_n\| + \|A_{n+1}x_{n+1} - A_nx_{n+1}\|.
\end{align*}
\item We have that:
\begin{align*}
\|x_{n+1}-A_{n+1}x_{n+1}\| &= t_n\|x_n-A_{n+1}x_{n+1}\| + (1-t_n)\|A_nx_n - A_{n+1}x_{n+1}\|\\
&\leq t_n(1-t_n)\|x_n-A_nx_n\| + t_n\|x_{n+1}-A_{n+1}x_{n+1}\| + \\
&\qquad +(1-t_n)^2\|x_n-A_nx_n\| + (1-t_n)\|A_{n+1}x_{n+1} - A_nx_{n+1}\| \\
&\leq (1-t_n)\|x_n-A_nx_n\| + t_n\|x_{n+1}-A_{n+1}x_{n+1}\|  +\\
&\qquad + (1-t_n)\|A_{n+1}x_{n+1} - A_nx_{n+1}\|.
\end{align*}
\item Immediate from the last statement.
\ee
\end{proof}

\bthm[``rate of $A_n$-asymptotic regularity'']\label{Than}
Set, for all $\eps>0$:
\begin{align*}
\Phi_{b,\theta,\gamma}(\eps):=&\ \Delta_{b,\theta}\left(\frac{\eps}{2},\gamma\left(\frac{\eps}{6b}\right)+1\right)\\
=&\ \theta\left(\gamma\left(\frac{\eps}{6b}\right)+\left\lceil\frac{4b^2}{\eps^2}\right\rceil+1\right).
\end{align*}
Then for all $\eps>0$ we have that for all $n \geq \Phi_{b,\theta,\gamma}(\eps)$, $\|x_n-A_nx_n\| \leq \eps$.
\ethm

\begin{proof}
First we get, for any $i$ and $n$, that:
\begin{align*}
\|T_ix_{n+1}\|&\leq \|T_ix_{n+1}-p\|+\|p-x\|+\|x\| \\
&\leq \|x_{n+1}-p\|+\|p-x\|+\|x\| \\
&\leq \|x-p\|+\|p-x\|+\|x\| \\
&\leq 3b.
\end{align*}
Using Lemma~\ref{lc}, we get that:
\begin{align*}
\|x_{n+1} - A_{n+1}x_{n+1}\| &\leq \|x_n-A_nx_n\| + \|A_{n+1}x_{n+1} - A_nx_{n+1}\| \\
&\leq \|x_n-A_nx_n\| + \|\sum_{i=1}^N (\lambda_i^{(n+1)}-\lambda_i^{(n)}) T_ix_{n+1} \| \\
&\leq \|x_n-A_nx_n\| + 3b \sum_{i=1}^N |\lambda_i^{(n+1)}-\lambda_i^{(n)}|.
\end{align*}
and, by summing up, we have that for any $n,m \in \N$:
\beq\label{sum}
\|x_{n+m} - A_{n+m}x_{n+m}\| \leq \|x_n-A_nx_n\| + 3b\sum_{j=n}^{n+m-1} \sum_{i=1}^N |\lambda_i^{(j+1)}-\lambda_i^{(j)}|.
\eeq
We now apply Proposition~\ref{liminf} for $\frac{\eps}{2}$ and $\gamma(\frac{\eps}{6b})+1$ to get that there is an $N$ in the interval
$$[\gamma(\frac{\eps}{6b})+1,\Delta_{b,\theta}(\frac{\eps}{2},\gamma(\frac{\eps}{6b})+1)=\Phi_{b,\theta,\gamma}(\eps)]$$ such that:
$$\|x_N - A_Nx_N\|\leq\frac{\eps}{2}.$$
Take now an arbitrary $n \geq \Phi_{b,\theta,\gamma}(\eps)$. Then, applying \eqref{sum}, we have that:
\begin{align*}
\|x_n - A_nx_n\|&=\|x_{N+(n-N)} - A_{N+(n-N)}x_{N+(n-N)}\|\\
&\leq \|x_N - A_Nx_N\| + 3b\sum_{j=N}^{n-1} \sum_{i=1}^N |\lambda_i^{(j+1)}-\lambda_i^{(j)}| \\
&\leq \|x_N - A_Nx_N\| + 3b\sum_{j=\gamma(\frac{\eps}{6b})+1}^{n-1} \sum_{i=1}^N |\lambda_i^{(j+1)}-\lambda_i^{(j)}| \\
&\leq \|x_N - A_Nx_N\| + 3b \cdot \frac{\eps}{6b}\\
&\leq \frac{\eps}{2} + \frac{\eps}{2} = \eps,
\end{align*}
where at the penultimate line we used the fact that $\gamma$ is a Cauchy modulus. The statement is now proven.
\end{proof}

Note that using this alternate proof we have eliminated the square root from the original hypotheses of López-Acedo and Xu. This will be maintained when using the statement above for $k$-strict pseudocontractions, and also for the general convergence theorem, which does not use that specific hypothesis in its proof except in the preliminary asymptotic regularity result.

\subsection{Individual rates of asymptotic regularity}

As announced earlier, we shall derive the rates of asymptotic regularity for each of the $T_i$'s.

Define the function $h : (0,\infty) \to (0,\infty)$, for any $\eps>0$, as follows:
$$h_{a,b}(\eps):=\eps + \sqrt{\frac{1-a}{a}}\sqrt{\eps^2+2b\eps}$$

\blem\label{preti}
Let $n \in \N$ and $z \in C$ such that $\|z-p\| \leq b$ and $\|z - A_nz\| \leq \eps$. Then, for each $i$, $\|z-T_iz\| \leq h_{a,b}(\eps)$.
\elem

\begin{proof}
Remember that $A_n = \sum_{i=1}^N \lambda_i^{(n)} T_i$ and that $p$ is a fixed point of all the $T_i$'s and hence also of the $A_n$'s. Then:
$$\|p-A_nz\| = \|A_np - A_nz\| \leq \|p-z\| \leq b.$$
Let $j \in \{1,...,N\}$. We have that:
\begin{align*}
\|p-z\|^2 &\leq \|p-A_nz\|^2 + \|z-A_nz\|^2 + 2\|p-A_nz\|\|z-A_nz\| \\
 &\leq \left\|p-\left(\left(\sum_{i \neq j} \lambda_i^{(n)}\right)\left(\frac{1}{\sum_{i \neq j} \lambda_i^{(n)}}\sum_{i \neq j} \lambda_i^{(n)} T_iz\right) + \lambda_j^{(n)}T_jz\right)\right\|^2 + \|z-A_nz\|^2 + \\
 &\qquad + 2\|p-A_nz\|\|z-A_nz\| \\
 &\leq \left(\sum_{i \neq j} \lambda_i^{(n)}\right)\left\|p - \frac{1}{\sum_{i \neq j} \lambda_i^{(n)}}\sum_{i \neq j} \lambda_i^{(n)} T_iz \right\|^2 + \lambda_j^{(n)}\|p-T_jz\|^2 -\\
 &\qquad - \lambda_j^{(n)}\left(\sum_{i \neq j} \lambda_i^{(n)}\right)\left\|\frac{1}{\sum_{i \neq j} \lambda_i^{(n)}}\sum_{i \neq j} \lambda_i^{(n)} T_iz - T_jz\right\|^2 + \eps^2 +2b\eps \\
 &\leq \left(\sum_{i \neq j} \lambda_i^{(n)}\right)\left\|p - z\right\|^2 + \lambda_j^{(n)}\|p-z\|^2 - \\
 &\qquad - \lambda_j^{(n)}\left(\sum_{i \neq j} \lambda_i^{(n)}\right)\left\|\frac{1}{\sum_{i \neq j} \lambda_i^{(n)}}\sum_{i \neq j} \lambda_i^{(n)} T_iz - T_jz\right\|^2 + \eps^2 +2b\eps \\
 &\leq \|p-z\|^2 - \lambda_j^{(n)}\left(\sum_{i \neq j} \lambda_i^{(n)}\right)\left\|\frac{1}{\sum_{i \neq j} \lambda_i^{(n)}}\sum_{i \neq j} \lambda_i^{(n)} T_iz - T_jz\right\|^2 + \eps^2 +2b\eps. \\
\end{align*}
It follows that:
$$\lambda_j^{(n)}\left(\sum_{i \neq j} \lambda_i^{(n)}\right)\left\|\frac{1}{\sum_{i \neq j} \lambda_i^{(n)}}\sum_{i \neq j} \lambda_i^{(n)} T_iz - T_jz\right\|^2 \leq \eps^2 +2b\eps,$$
so
$$\left\|\frac{1}{\sum_{i \neq j} \lambda_i^{(n)}}\sum_{i \neq j} \lambda_i^{(n)} T_iz - T_jz\right\| \leq \sqrt{\frac{1}{\lambda_j^{(n)}\left(\sum_{i \neq j} \lambda_i^{(n)}\right)}}\sqrt{\eps^2 +2b\eps}.$$
We can now see that:
\begin{align*}
\|z-T_jz\| &\leq \|z-A_nz\| + \|A_nz - T_jz\| \\
&\leq \eps + \left(\sum_{i \neq j} \lambda_i^{(n)}\right)\left\|\frac{1}{\sum_{i \neq j} \lambda_i^{(n)}}\sum_{i \neq j} \lambda_i^{(n)} T_iz - T_jz\right\| \\
&\leq \eps+\sqrt{\frac{\sum_{i \neq j} \lambda_i^{(n)}}{\lambda_j^{(n)}}}\sqrt{\eps^2 +2b\eps}.
\end{align*}
All that remains is to show that
$$\frac{\sum_{i \neq j} \lambda_i^{(n)}}{\lambda_j^{(n)}} \leq \frac{1-a}{a}.$$
But this is immediate, since $\lambda_j^{(n)} \geq a$, and hence $\frac{1}{\lambda_j^{(n)}} \leq \frac{1}{a}$ and $\sum_{i \neq j} \lambda_i^{(n)} = 1- \lambda_j^{(n)} \leq 1-a$.
\end{proof}

\bthm[``rate of $T_i$-asymptotic regularity'']\label{Tinonexp}
Set, for all $\eps>0$:
\begin{align*}
P_{a,b}(\eps):=& \min \left\{ \frac{\eps}{2} , \sqrt{\frac{a\eps^2}{4(1-a)}+b^2} -b \right\} \\
\Phi'_{a,b,\theta,\gamma}(\eps):=&\ \Phi_{b,\theta,\gamma}(P_{a,b}(\eps)).
\end{align*}
Then for all $\eps>0$ we have that for all $n \geq \Phi'_{a,b,\theta,\gamma}(\eps)$ and all $i$, $\|x_n-T_ix_n\| \leq \eps$.
\ethm

\begin{proof}
By Theorem~\ref{Than}, we have that $\|x_n - A_nx_n\| \leq P_{a,b}(\eps)$. Since $\|x_n - p \| \leq \|x-p\| \leq b$, we can apply Lemma~\ref{preti} to get that:
$$\|x_n - T_ix_n\| \leq h_{a,b}(P_{a,b}(\eps)).$$
All we have to show, therefore, is that $h_{a,b}(P_{a,b}(\eps)) \leq \eps$, i.e. that:
$$ P_{a,b}(\eps) + \sqrt{\frac{1-a}{a}}\sqrt{ P_{a,b}(\eps)^2+2b P_{a,b}(\eps)} \leq \eps.$$
Now, it is immediate that the monic polynomial:
$$X^2 + 2bX - \frac{a\eps^2}{4(1-a)}$$
has two roots, one negative and one equal to $\sqrt{\frac{a\eps^2}{4(1-a)}+b^2} -b$. Since (by its definition) $P_{a,b}(\eps) \geq 0$ and $P_{a,b}(\eps) \leq \sqrt{\frac{a\eps^2}{4(1-a)}+b^2} -b$, we have that:
$$P_{a,b}(\eps)^2 + 2bP_{a,b}(\eps) - \frac{a\eps^2}{4(1-a)} \leq 0,$$
and hence that
$$P_{a,b}(\eps)^2 + 2bP_{a,b}(\eps) \leq \frac{a\eps^2}{4(1-a)}.$$
We can now compute:
\begin{align*}
P_{a,b}(\eps) + \sqrt{\frac{1-a}{a}}\sqrt{ P_{a,b}(\eps)^2+2b P_{a,b}(\eps)} &\leq \frac{\eps}{2} + \sqrt{\frac{1-a}{a}}\sqrt{\frac{a\eps^2}{4(1-a)}} \\
&= \frac{\eps}{2} + \frac{\eps}{2} = \eps.
\end{align*}
\end{proof}

\subsection{The general case}

We can at last consider the general case of a family of strict pseudocontractions -- i.e. we will allow $k$ to vary. Note that if $0 \leq k \leq k' <1$ and $T$ is a $k$-strict pseudocontraction, then it is also a $k'$-strict pseudocontraction. Hence, instead of considering each $T_i$ to be a $k_i$-strict pseudocontraction, we can take $k$ to be the maximum of the $N$ $k_i$'s and work, without loss of generality, with a finite family of $k$-strict pseudocontractions.

Now, if $T : C \to C$ is a $k$-strict pseudocontraction and $t \in (0,1)$, denote by $T_t : C \to C$ the map defined by $T_t x := tx + (1-t)x$. It is clear that for any such $t$, the maps $T$ and $T_t$ have the same fixed points. We shall make use of the following lemma. (See \cite{Sip16} for more uses of this kind of argument to establish quantitative and qualitative results for strict pseudocontractions.)

\blem[{\cite[Theorem 2]{BroPet67}}]\label{bpl}
If $T : C \to C$ is a $k$-strict pseudocontraction, then $T_k$ is nonexpansive.
\elem

Also, note that $(T_{t_1})_{t_2} = T_{1- (1-t_1)(1-t_2)}$. (We have used here, as in the iterative algorithms, the ``opposite'' convention from the one in \cite{BroPet67}, namely the one in \cite{MarXu07,LopXu07}, hence the form of this result.)

Our main result will be the following.

\bthm[``general rate of $T_i$-asymptotic regularity'']
Set, for all $\eps>0$:
\begin{align*}
\Phi''_{a,b,k,\theta,\gamma}(\eps):=&\ \Phi'_{a,b,\theta,\gamma}((1-k)\eps).
\end{align*}
Then for all $\eps>0$ we have that for all $n \geq \Phi''_{a,b,k,\theta,\gamma}(\eps)$ and all $i$, $\|x_n-T_ix_n\| \leq \eps$.
\ethm

\begin{proof}
Set, for all $i$, $T'_i:=(T_i)_k$. Then, by Lemma~\ref{bpl}, each such $T'_i$ is nonexpansive. It is easy to see that if we set $A'_n := \sum_{i=1}^N \lambda_i^{(n)} T'_i$, then $A'_n = (A_n)_k$.

Set, for all $n$, $t'_n := 1-\frac{1-t_n}{1-k}$. Since $t_n \in [k,1]$, we have that $t'_n \in [0,1]$.

Let $n \in \N$. Then we have that:
\begin{align*}
x_{n+1} &= t_nx_n + (1-t_n)A_nx_n \\
&= (A_n)_{t_n} x_n \\
&= (A_n)_{1-(1-k)(1-t'_n)} x_n \\
&= ((A_n)_k)_{t'_n} x_n \\
&= (A'_n)_{t'_n} x_n.
\end{align*}
So the sequence $(x_n)_n$ is also the output of the parallel algorithm for the mappings $(T'_i)_i$ and the weight sequence $(t'_n)_n$, with the same $\lambda_i^{(n)}$'s. Now, in order to apply Theorem~\ref{Tinonexp}, we must check that the $\theta$ is still valid. We have that $1-t_n = (1-k)(1-t'_n)$ and that $t_n-k=(1-k)t'_n$, so:
$$\sum_{n=0}^{\theta(N)} t'_n(1-t'_n) = \frac{1}{(1-k)^2} \sum_{n=0}^{\theta(N)} (t_n -k)(1-t_n) \geq \frac{1}{(1-k)^2} \cdot N \geq N.$$
Hence, by Theorem~\ref{Tinonexp} we have that for all $\eps>0$, for all $n \geq \Phi'_{a,b,\theta,\gamma}(\eps)$ and for all $i$, $\|x_n-T'_ix_n\| \leq \eps$. But $x_n -T'_i x_n = (1-k)(x_n -T_i x_n)$, so:
$$\|x_n -T_i x_n \| \leq \frac{\eps}{1-k}.$$

So, for all $n \geq \Phi'_{a,b,\theta,\gamma}((1-k)\eps)$, we have that:
$$\|x_n -T_i x_n \| \leq \eps,$$
which was what we needed to show.
\end{proof}

As a final aside, we note that one can remove the explicit appearance of $k$ in the above bound by the following argument. By the definition of $\theta$, we have that:
$$\sum_{n=0}^{\theta(1)} (t_n-k)(1-t_n) \geq 1.$$
Therefore, by the pigeonhole principle, we have that one of these first $\theta(1)+1$ terms of the series must be greater than $\frac{1}{\theta(1)+1}$. So there is an $n_0 \leq \theta(1)$ such that:
$$(t_{n_0}-k)(1-t_{n_0}) \geq \frac{1}{\theta(1)+1}.$$
But necessarily $(t_{n_0}-k)(1-t_{n_0}) \leq t_{n_0}-k \leq 1-k$, so $k \leq 1 - \frac{1}{\theta(1)+1}$. We noted earlier that we can replace any $k$ by a greater $k'$, smaller than $1$, and the strict pseudocontraction condition will be maintained. Therefore, in the above bound, we may freely replace $k$ by $1 - \frac{1}{\theta(1)+1}$ (taking note, in the process, that $\theta(1)$ plays the r\^ole of $M$ from the statements of the two metatheorems).

\section{Acknowledgments}

The author is grateful to Lauren\c tiu Leu\c stean for originally suggesting the problem and to him and Ulrich Kohlenbach for their comments that greatly contributed to the final form of this paper.

This work was supported by a grant of the Romanian National Authority for Scientific Research, CNCS - UEFISCDI, project number PN-II-ID-PCE-2011-3-0383.

\end{document}